\documentclass[11pt,reqno]{amsart}

\usepackage{latexsym}
\usepackage{amsfonts}
\usepackage{amsmath}
\usepackage{url}
\usepackage{graphicx}
\usepackage{color}

\usepackage{epsfig,psfrag,verbatim,amsfonts}

\usepackage{amssymb,amsmath}

\usepackage{xcolor}
\usepackage[first=375, last=475]{lcg}

\def\bl{\begin{lemma}}
\def\el{\end{lemma}}
\def\bth{\begin{theorem}}
\def\eth{\end{theorem}}
\def\bc{\begin{corollary}}
\def\ec{\end{corollary}}
\def\bcj{\begin{conjecture}}
\def\ecj{\end{conjecture}}
\def\bpr{\begin{proposition}}
\def\epr{\end{proposition}}
\def\bde{\begin{definition}}
\def\ede{\end{definition}}

\newcommand{\be}{\begin{eqnarray}}
\newcommand{\ee}{\end{eqnarray}}

\renewcommand{\and}{\hbox{ {\rm and} }}

\newtheorem{theorem}{Theorem}[section]
\newtheorem{definition}{Definition}[section]
\newtheorem{lemma}[theorem]{Lemma}

\newtheorem{corollary}[theorem]{Corollary}
\newtheorem{proposition}[theorem]{Proposition}
\newtheorem{conjecture}[theorem]{Conjecture}

\theoremstyle{definition}
\numberwithin{equation}{section}

\begin{document}
\title{Sensitivity of mixing times, an example}
\author{Itai Benjamini  }

\begin{abstract}
We construct a family of growing finite bounded degree rooted graphs, $G_n$,
in which the mixing time for simple random walk, starting at the root, is order
$\log |G_n|$. Yet after a  quasi - isometry, the ratio of $|G_n|$ over the mixing time
grows  arbitrarily slow to infinity.

\end{abstract}

\maketitle

\section{Introduction}

The question of instability of the mixing time was thoroughly studied in \cite{DP, H, HK, HP}. Sharp bounds were established on the ratio of the mixing times between two quasi-isometric graphs, for several notions of mixing times, total variation distance, $l^2$ and $l^{\infty}$, when the walk started from the worse starting vertex. The spectral gap is quasi - isometric invariant. By the relation between mixing times and spectral gap,  when the degree is bounded, the ratio of the mixing times when starting from worse starting vertex can be at most $\log$ the size of the graph.

In this note we construct a growing family of rooted bounded degree  graphs, denoted $G_n$,  were by a quasi-isometry any of the above notions of mixing times, jumps from
$ \log |G_n|$ to $|G_n| \over{ f(n) }$, for a function $f(n)$ that grows arbitrarily slow to infinity.
Such a gap can not be obtained in a vertex transitive graph. Let me advertise the remarkable, yet involved vertex transitive  example in \cite{HK}.
Is this  the largest possible  ratio? that can be achieved? 

\medskip

We now recall the definition of the total variation mixing time and quasi - isometry.

For any two distributions $\mu$ and $\nu$  on a finite state space, e.g. the vertices of the graph, their total variation distance is defined as
$1/2 \sum_{x} |\mu(x) - \nu(x)|$. The total variation distance for a random walk starting at $x$ at time $t$ is defined
as $d^x(t)  =  ||P^t_x - \pi||_{TV}$. Where by $P^t_x$ we denote
the distribution of the simple random walk starting at $x$ by time $t$ and $\pi$ is the stationary distribution.
The $\epsilon$ - total variation mixing time for a random walk starting at $x$ is defined as
$ t_{mix}(\epsilon)  = \inf \{t : d^x(t) < \epsilon\}$.
\medskip

A {\it quasi-isometry} between
two metric spaces $X$ and $Y$ is a map $\psi : X \rightarrow Y$ such that for some numbers $(a, b)$
we have $\forall u, v \in  X$

$a^{-1} d(u,v) - b \leq d(\psi(u), \psi(v)) \leq a d(u,v) + b$.

where $d$ denotes the distance (in $X$ or in $Y$ , as appropriate); and further such that
for every $y \in  Y$ there is some $x \in X$ such that $d(\psi(x), y) \leq a + b$. We say that
$X$ and $Y$ are $(a, b)$ quasi-isometric if such a $\psi$ exists.

In our setup, we will have fixed $(a,b)$ and a sequence of pairs
of bounded degree graphs, growing in size, in which the graphs of  each pair are $(a,b)$ quasi isometric.

\section{The example}
A fact used by all the constructions cited above (\cite{DP, H, HK, HP}) is that in the  rooted infinite binary tree (or any $d > 1$ - regular tree, there is a subset of vertices, deonted by $A$, such that simple random walk on the binary tree will visit infinitely often with probability one. Yet after a quasi - isometry of the tree, the random walk on the modified tree starting from the root  will not hit $A$ with probability that can be made arbitrarily close to one. Denote by $A_n$ the intersection of $A$ with the ball of radius $n$ around the root. To get $A$,  take all the vertices so that in their representation as words in $0,1$, for turning left or right along the path from the root, when drawing the tree in the plane, the ratio of the number of $0$ to the number of $1$ is between $0.9999$ and $1.0001$. For the  quasi-isometry replace each right turning edge by a path of length ten. Such a quasi - isometry will asymptotically change the density of $1$'s in the vertices visited by the simple random walk.
\medskip

Consider a family of growing $d$-regular, vertex transitive, large girth expanders, denote by $H_n$. Were $H_n$ is chosen such that balls of radius $\Omega(n)$ (say $2^n$) are trees. E.g. Ramanujan graphs \cite{HLW}.

To construct the example, denoted $G_n$, for each $n$, take two copies of $H_n$ pick a root in each of the copies. A ball of radius $n$ centred at each of the roots, in the respective copies is a $d$-tree up until the leaves. By isomorphically mapping the ball to the ball in the  infinite tree  choose an identification of  $A_n$ in each of the copies.  To get $G_n$ connect, {\it by a path of length two}, each vertex in the identification of $A_n$ in one of the copies to it's identical vertex in the identification of $A_n$ in the other copies. Note that the ball around the root that contained $A_n$ has size which is $o(|G_n|)$. This is the example, see figure.

The quasi - isometry will be replacing in the $n$  - ball around the roots,  each right turning edge by a path of length ten.
\medskip

\begin{figure}
  \includegraphics[width=\linewidth]{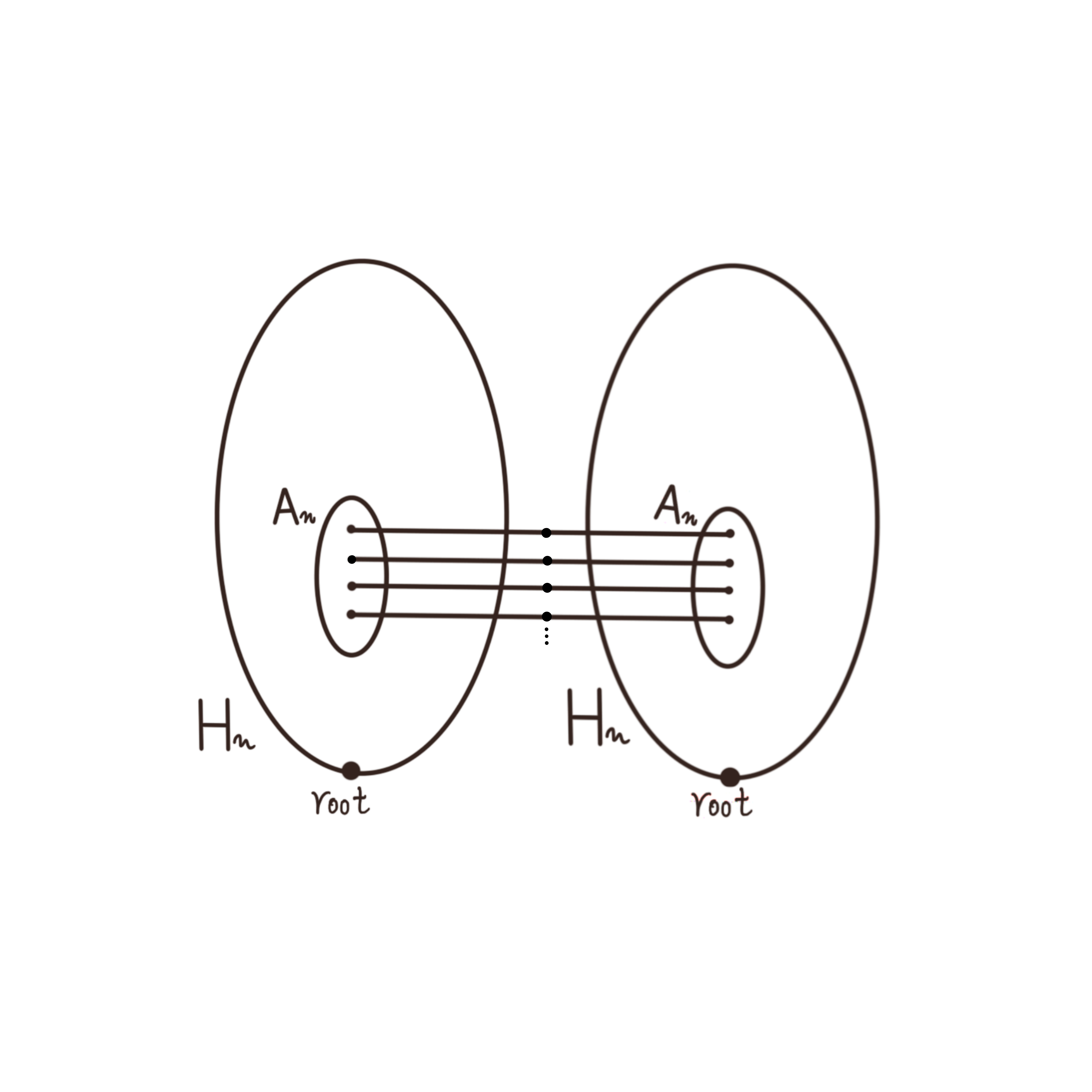}
  \caption{The example.}
\end{figure}

Note that for any $k >0$, with probability going to one with $n$, simple random walk starting from the root will visit  $A_n$ is each of the copies more than $k$ times.  Thus with high probability the simple random walk will cross to the other copy. Couple two simple random walks, one starting from each of the roots, so that the walks takes identical steps, each walker in her copy, and move to the middle vertex on a length two paths between elements of $A_n$ together.  Once they meet in the middle of one of the length two paths, let them continue together.  This event  happens  with probability arbitrarily close to one before reaching distance $n$ from the respective roots. Conditioning on meeting, the probability the walkers will absorbed after exiting the $n$-balls in each of the copies, is $1/2$. Therefore the mixing time of $G_n$, conditioning on meeting is the mixing time, is that of the expander $H_n$, which is $\log |H_n|$.
\medskip

After applying the quasi - isometry in the  ball of radius $n$ around the roots, with probability which can go to one with $n$, Simple random walk starting at any of the roots, will exit the $n$-ball before hitting $A_n$ and will mix in the copy it started. Simple random walk on an expander mixes in order  $\log$ size steps, which happens with high probability, before returning to a polynomially small fixed ball if starting at distance growing to infinity with $n$ from the ball. Therefore the hitting time of $A_n$ from a uniform starting vertex in the copy is a lower bound on the mixing time of $G_n$, which is bounded from below by $|G_n| \over{|A_n|}$.
To sse this let $T$ be the hitting time of the set $S$
let $V_t(S)$ be the number of visits to $S$ up to time $t$.
Starting from stationary distribution we have by Markov's inequality,
$$
 \mathbb{P} [ T \leq t ] \leq \mathbb{P} [ V_t(S) \geq 1 ] \leq \mathbb{E} [ V_t(S) ] = \pi(S) t
$$
where $\pi$ is the stationary measure
if $\pi$ is uniform (for regular graphs), then this will be small as long as $t < \frac{1}{\pi(S)} = \frac{N}{|S|}$
\qed
\medskip

The properties of  $H_n$ we used are that as $n$ grow, balls of radius growing to infinity around the roots are trees, and that $H_n$ has logarithmic mixing time.

\noindent{\bf Remark:} As we commented  a key to all the examples is the fact that in the infinite binary tree, there is a set that absorbeds the simple random walk. Yet after a quasi - isometry simple random a.s. visits the set finitely many times. In \cite{BR} it was shown that such a set
exists on some amenable groups including  the lamplighter over $Z$, with respect to two symmetric generating sets. It is of interest to know what other Cayley graphs admits such an unstable absorbing set? Does all super polynomial Cayley graphs admit an unstable absorbing set?

\noindent
{\bf Acknowledgements:} Thanks to Jonathan Hermon, Elad Tzalik and  Ariel Yadin for useful discussions.

\end{document}